\documentclass{article}

\usepackage[applemac]{inputenc} 
\usepackage[dvips]{graphicx}
\usepackage{amsmath} 
\usepackage{amssymb} 

\makeatletter
\def\@maketitle{%
  \vbox to 3.0in{%
    \hsize\textwidth
    \linewidth\hsize
    \vspace*{1.5cm}
    \centering
    {\bfseries\huge \@title \par}
    \vskip 2em
    {\large \begin{tabular}[t]{c}\@author \end{tabular}\par}
    \vfill}    \vspace*{1.0cm}
}
\renewcommand\section{\@startsection {section}{1}{\z@}%
     {.7\baselineskip plus\baselineskip}{.5\baselineskip}
                                   {\normalfont\Large\bfseries}}
\renewcommand\section{\@startsection {section}{1}{\z@}%
      {.5\baselineskip\@plus.7\baselineskip}{.3\baselineskip}%
                                   {\normalfont\Large\bfseries}}
\renewcommand\subsection{\@startsection{subsection}{2}{\z@}%
       {.5\baselineskip\@plus.7\baselineskip}{.3\baselineskip}%
                                   {\normalfont\large\bfseries}}
\renewcommand\subsubsection{\@startsection{subsubsection}{3}{\z@}%
      {.5\baselineskip\@plus.7\baselineskip}{.3\baselineskip}%
                                     {\normalfont\normalsize\bfseries}}
\makeatother

\renewenvironment{abstract}%
  {\normalfont
    \list{}{\labelwidth0pt
      \leftmargin0pt \rightmargin\leftmargin
      \listparindent\parindent \itemindent0pt
      \parsep0pt
      
    }%
    \item[\hskip\labelsep\bfseries\abstractname\enspace --] \itshape%
}{%
  \endlist}

\newcommand{\keywordsname}{Keywords}
\newenvironment{keywords}%
  {\normalfont
    \list{}{\labelwidth0pt
      \leftmargin0pt \rightmargin\leftmargin
      \listparindent\parindent \itemindent0pt
      \parsep0pt
      }%
    \item[\hskip\labelsep\bfseries\keywordsname:]}{\endlist}

\begin{document}

\pagestyle{plain} 

\title{On the Blackman's Association Problem}

\author{Jean Dezert\\
ONERA\\
29 Avenue de la  Division Leclerc \\
92320 Ch\^{a}tillon, France.\\
Jean.Dezert@onera.fr\\
\and
Florentin Smarandache\\
Department of Mathematics\\
University of New Mexico\\
Gallup, NM 87301, U.S.A.\\
smarand@unm.edu
\and
Albena Tchamova\\
CLPP, Bulgarian Academy of Sciences\\
Sofia, Bulgaria.\\
tchamova@bas.bg}

\date{}

\maketitle

\begin{abstract}
Modern multitarget-multisensor tracking systems involve the development of reliable methods for the data association and the fusion of multiple sensor information, and more specifically the partioning of observations into tracks. This paper discusses and compares the application of Dempster-Shafer Theory (DST) and the Dezert-Smarandache Theory (DSmT) methods to the fusion of multiple sensor attributes for target identification purpose. We focus our attention on the paradoxical Blackman's association problem and propose several approaches to outperfom Blackman's solution. We  clarify some preconceived ideas about the use of degree of conflict between sources as potential criterion for partitioning  evidences.
\end{abstract}

\begin{keywords}
Data Association, Entropy, Data Fusion, Uncertainty, Paradox, Dezert-Smarandache theory, plausible and paradoxical reasoning
\end{keywords}

\section{Introduction}

The association problem is of major importance in most of modern multitarget-multisensor tracking systems. This task is particularly difficult when data are uncertain and are modeled by basic belief masses and when sources are conflicting. The solution adopted is usually based on the Dempster-Shafer Theory (DST)  \cite{Shafer_1976} because it provides an elegant theoretical way to combine uncertain information. However the Dempster's rule of combination can give rise to some paradox/anomaly and can fail to provide the correct solution for some specific association problems. This has been already pointed out  by Samuel Blackman in \cite{Blackman_1990}. Therefore more study in this area is required and we propose here a new analysis of the Blackman's association problem (BAP). We present in the sequel the original BAP and remind the classical attempts to solve it based on DST (including the Blackman's method). In the second part of the paper we propose and compare new approches based on the recent Dezert-Smarandache Theory (DSmT) of plausible and paradoxical reasoning \cite{Dezert_2002b,Smarandache_2002}. The DSmT can be interpreted as a generalization of the DST and allows to combine formally any types of sources of information (rational, uncertain or paradoxical). The last part of the paper provides a comparison of the performances of all the proposed approaches from Monte-Carlo simulation results.

\section{The Association Problem}

\subsection{Association Problem no. 1}

Let's recall now the original Blackman's association problem \cite{Blackman_1990}. Consider only two target attribute types 
corresponding to the very simple {\sl{frame of discernment}}
$\Theta=\{\theta_{1},\theta_{2}\}$ and the association/assignment problem for a 
single attribute observation $Z$ and two tracks ($T_{1}$ and $T_{2}$). Assume now 
the following two predicted basic belief assignments (bba) for attributes of the two tracks:
$$m_{T_{1}}(\theta_{1})=0.5\quad m_{T_{1}}(\theta_{2})=0.5\quad m_{T_{1}}(\theta_{1}\cup\theta_{2})=0$$
$$m_{T_{2}}(\theta_{1})=0.1\quad m_{T_{2}}(\theta_{2})=0.1\quad m_{T_{2}}(\theta_{1}\cup\theta_{2})=0.8$$
We now assume to receive the new following bba drawn from attribute 
observation $Z$ of the system
$$m_{Z}(\theta_{1})=0.5\quad m_{Z}(\theta_{2})=0.5\quad m_{Z}(\theta_{1}\cup\theta_{2})=0$$
The problem is to develop a general method to find the correct assignment of the attribute measure $m_{Z}(.)$ with the predicted one $m_{T_{i}}(.)$, $i=1,2$. 
Since $m_{Z}(.)$ matches perfectly with $m_{T_{1}}(.)$ 
whereas $m_{Z}(.)$ does not match with $m_{T_{2}}(.)$, the optimal solution is obviously given by the assignment $(m_{Z}(.) \leftrightarrow m_{T_{1}}(.))$. The problem is to find an unique general and reliable method for solving this specific problem and for solving all the other possible association problems as well.

\subsection{Association Problem no. 2}

To compare several potential issues, we propose to modify the previous problem into a second one by keeping the same predicted bba $m_{T_{1}}(.)$ and $m_{T_{2}}(.)$ but by considering now the following bba $m_{Z}(.)$
$$m_{Z}(\theta_{1})=0.1\quad m_{Z}(\theta_{2})=0.1\quad m_{Z}(\theta_{1}\cup\theta_{2})=0.8$$
Since $m_{Z}(.)$ matches perfectly with $m_{T_{2}}(.)$, the correct solution is now directly given by 
$(m_{Z}(.) \leftrightarrow m_{T_{2}}(.))$. The sequel of this paper in devoted to the presentation of some attempts for solving the BAP, not only for these two specific problems 1 and 2, but for the more general problem where the bba $m_Z(.)$ does not match perfectly with one of the predicted bba $m_{T_i}$, $i=1$ or $i=2$ due to observation noises.

\section{Attempts for solutions}

We examine now several approaches which have already been (or could be) envisaged 
 to solve the general association problem.

\subsection{The simplest approach}

The simplest idea for solving BAP, surprisingly not reported by Blackman in \cite{Blackman_1990} is to use a classical minimum distance criterion directly between the predictions $m_{T_i}$ and the observation $m_Z$. The classical $L^1$ (city-block) or $L^2$ (Euclidean) distances are typically used. Such simple criterion obviously provides the correct association in most of cases involving perfect (noise-free) observations $m_Z(.)$.  But there exists numerical cases for which the optimal decision cannot be found at all, like in the following numerical example:
$$m_{T_{1}}(\theta_{1})=0.4\quad m_{T_{1}}(\theta_{2})=0.4\quad m_{T_{1}}(\theta_{1}\cup\theta_{2})=0.2$$
$$m_{T_{2}}(\theta_{1})=0.2\quad m_{T_{2}}(\theta_{2})=0.2\quad m_{T_{2}}(\theta_{1}\cup\theta_{2})=0.6$$
$$m_{Z}(\theta_{1})=0.3\quad m_{Z}(\theta_{2})=0.3\quad m_{Z}(\theta_{1}\cup\theta_{2})=0.4$$
\noindent
From these bba, one gets $d_{L^1}(T_{1},Z) = d_{L^1}(T_{2},Z)=0.4$ (or $d_{L^2}(T_{1},Z) = d_{L^2}(T_{2},Z)\approx0.24$) and no decision can be drawn for sure,  although the minimum conflict approach (detailed in next section) will give us instead the following solution $(Z\leftrightarrow T_2)$. It is not obvious in such cases to justify this method with respect to some other ones.
What is more important in practice \cite{Blackman_1990}, is not only the association solution itself but also the attribute likelihood function $P(Z|T_i)\equiv P(Z\leftrightarrow T_i)$. As we know many likelihood functions (exponential, hyperexponential, Chi-square, Weibull pdf, etc) could be build from $d_{L^1}(T_{i},Z)$ (or $d_{L^2}(T_{i},Z)$ measures but we do not  know in general which one corresponds to the real attribute likelihood function.

\subsection{The minimum conflict approach}

The first  idea suggested by Blackman for solving the association problem was to apply the Dempster's rule of combination \cite{Shafer_1976} $m_{{T_i}Z}(.)=[m_{T_i}\oplus m_{Z}](.)$ defined by $m_{{T_i}Z}(\emptyset)=0$ and for any $C\neq\emptyset$ 
 and $C\subseteq\Theta$,
 \begin{equation*}
 m_{{T_i}Z}(C)= \frac{1}{1-k_{T_iZ}} \sum_{A\cap B=C}m_{T_i}(A)m_{Z}(B)
 \end{equation*}
and choose the solution corresponding to the minimum of conflict $k_{T_iZ}$. The sum in previous formula is over all $A,B\subseteq\Theta$ such that $A\cap B=C$. The degree of conflict $k_{T_iZ}$ between $m_{T_i}$ and $m_{Z}$ is given by $\sum_{A\cap B=\emptyset} m_{T_i}(A) m_{Z}(B)\neq 0$. Thus, an intuitive choice for the attribute likelihood function is $P(Z|T_i)=1-k_{T_iZ}$. If we now apply the Dempster's rule for the problem 1, we get the same result for both assignments, i.e.
$m_{T_{1}Z}(.)=m_{T_{2}Z}(.)$ with $m_{T_{i}Z}(\theta_{1})=m_{T_{i}Z}(\theta_{2})=0.5$ for $i=1,2$ and
$m_{T_{}Z}(\theta_{1}\cup\theta_{2})=0$, and more surprisingly, the correct assignment $(Z \leftrightarrow T_1)$ is not given by the minimum of conflict between sources since one has actually $(k_{T_{1}Z}=0.5)> (k_{T_{2}Z}=0.1)$. Thus, it is impossible to get the correct solution for this first BAP from the minimum conflict criterion as we firstly expected intuitively. This same criterion provides us however the correct solution for problem 2, since one has now $(k_{T_{2}Z}=0.02) < (k_{T_{1}Z}=0.1)$. The combined bba for problem 2 are given by $m_{T_{1}Z}(\theta_{1})=m_{T_{1}Z}(\theta_{2})=0.5$ and 
$m_{T_{2}Z}(\theta_{1})=m_{T_{2}Z}(\theta_{2})=0.17347$, $m_{T_{2}Z}(\theta_{1}\cup\theta_{2})= 0.65306$.

\subsection{The Blackman's approach}

To solve this apparent anomaly, Samuel Blackman has then proposed in \cite{Blackman_1990} to use a relative, 
rather than an absolute, attribute likelihood function as follows
$$L(Z\mid T_i)\triangleq (1-k_{T_iZ})/(1-k_{T_iZ}^{\text{min}})$$
\noindent
where $k_{T_iZ}^{\text{min}}$ is the minimum conflict factor that could 
occur for either the observation $Z$ or the track $T_i$ in the case of 
{\sl{perfect}} assignment (when $m_{Z}(.)$ and $m_{T_i}(.)$ coincide). By 
adopting this relative likelihood function, one gets now for problem 1

\begin{equation*}
\begin{cases}
L(Z\mid T_{1})=\frac{1-0.5}{1-0.5}=1\\
L(Z\mid T_{2})=\frac{1-0.1}{1-0.02}=0.92
\end{cases}
 \end{equation*}

Using this second Blackman's approach, there is now a larger likelihood 
associated with the first assignment (hence the right assignment 
solution for problem 1 can be obtained now based on the max likelihood criterion) but the difference between the 
two likelihood values is very small. As reported by S. 
Blackman in \cite{Blackman_1990}, {\sl{more study in this area is 
required}} and we examine now some other approaches. It is also interesting to note that  this same approach fails to solve the problem 2 since the corresponding likelihood functions 
for problem 2 become now
\begin{equation*}
\begin{cases}
L(Z\mid T_{1})=\frac{1-0.1}{1-0.5}=1.8\\
L(Z\mid T_{2})=\frac{1-0.02}{1-0.02}=1
\end{cases}
 \end{equation*}

\noindent
which means that the maximum likelihood solution gives now the incorrect assignment $(m_{Z}(.) \leftrightarrow m_{T_{1}}(.))$ for problem 2 as well, without mentioning the fact that the relative likelihood function becomes now greater than one !!!.

\subsection{The Tchamova's approach}

Following the idea of section 3.1, Albena Tchamova has recently proposed in \cite{Dezert_2003} to use rather  the $L^1$ (city-block) distance $d_{1}(T_i,T_iZ)$ or $L^2$ (Euclidean) distance $d_{2}(T_i,T_iZ)$ between the predicted bba $m_{T_i}(.)$ and the updated/combined bba $m_{T_iZ}(.)$ to measure the closeness of assignments with
$$d_{L^1}(T_i,T_iZ)=\sum_{A\in 2^\Theta}\mid m_{T_i}(A) - m_{T_iZ}(A)\mid$$
$$d_{L^2}(T_i,T_iZ)={[ \sum_{A\in 2^\Theta}[m_{T_i}(A) - m_{T_iZ}(A)]^2 ]}^{1/2}$$
The decision criterion here is again to choose the solution which yields  the minimum distance. This idea is  justified by the analogy with the steady-state Kalman filter (KF) behavior because if $z(k+1)$ and $\hat{z}(k+1| k)$ correspond to measurement and predicted measurement for time $k+1$, then the well-known KF updating state equation \cite{BarShalom_1993} is given by (assuming here that dynamic matrix is identity) $\hat{x}(k+1|k+1)=\hat{x}(k+1|k)+K(z(k+1)-\hat{z}(k+1| k))$. The steady-state is reached when $z(k+1)$ coincides with predicted measurement $\hat{z}(k+1| k)$ and therefore when $\hat{x}(k+1|k+1)\equiv\hat{x}(k+1|k)$. In our context, $m_{T_i(.)}$ plays the role of predicted state and $m_{T_iZ}(.)$ the role of updated state. Therefore it a priori makes sense that correct assignment should be obtained when $m_{T_iZ}(.)$ tends towards $m_{T_i}(.)$ for some closeness/distance criterion. Monte Carlo simulation results will prove however that this approach is also not as good as we can expect.\\

It is interesting to note that the Tchamova's approach succeeds to provide the correct solution for problem 1 with both distances criterions since $(d_{L^1}(T_{1},T_{1}Z) = 0) <   (d_{L^1}(T_{2},T_{2}Z) \sim 1.60)$ and $(d_{L^2}(T_{1},T_{1}Z) = 0) <   (d_{L^2}(T_{2},T_{2}Z) \sim 0.98)$, but provides the wrong solution for problem 2 since we will get both
$(d_{L^1}(T_{2},T_{2}Z)\sim 0.29) > (d_{L^1}(T_{1},T_{1}Z) = 0)$ and $(d_{L^2}(T_{2},T_{2}Z) \sim 0.18) > d_{L^2}(T_{1},T_{1}Z) = 0) $.

\subsection{The entropy approaches}

We examine here the results drawn from several entropy-like measures approaches. 
Our idea is now to use as decision criterion the minimum of the following entropy-like measures (expressed in {\it{nats}} - i.e. natural number basis with convention $0\log(0) = 0$):

\begin{itemize}
\item Extended entropy-like measure:
$$H_{ext}(m)\triangleq - \sum_{A\in2^\Theta} m(A)\log(m(A))$$
\item Generalized entropy-like measure \cite{Pal_1993,Schubert_2002}:
$$H_{gen}(m)\triangleq - \sum_{A\in2^\Theta} m(A)\log(m(A)/|A|)$$
\item Pignistic entropy:
$$H_{betP}(m)\triangleq - \sum_{\theta_i\in\Theta} P\{\theta_i\}\log(P\{\theta_i\})$$
\end{itemize}
\noindent where the pignistic(betting) probabilities $P(\theta_i)$ are obtained by $$\forall 
 \theta_{i}\in\Theta,  \quad P\{\theta_{i}\}=\sum_{B\subseteq\Theta|\theta_{i}\in B} 
 \frac{1}{|B|}m(B)$$

It can be easily verified that  the minimum entropy criterion (based on $H_{ext}$, $H_{gen}$ or $H_{betP}$) computed from combined bba 
$m_{{T1}Z}(.)$ or $m_{{T2}Z}(.)$ are actually unable to provide us correct solution for problem 1 because of indiscernibility of $m_{{T1}Z}(.)$ with respect to $m_{{T2}Z}(.)$. For problem 1, we get $H_{ext}(m_{{T1}Z})=H_{ext}(m_{{T2}Z})=0.69315$ and exactly same numerical results for $H_{gen}$ and $H_{betP}$ because no uncertainty is involved in the updated bba for this particular case. If we now examine the numerical results obtained for problem 2, we can see that minimum entropy
criteria is also unable to provide the correct solution based on $H_{ext}$, $H_{gen}$ or $H_{betP}$ criterions since 
one has $H_{ext}(m_{{T2}Z})= 0.88601> H_{ext}(m_{{T1}Z})=0.69315$, $H_{gen}(m_{{T2}Z})= 1.3387> H_{gen}(m_{{T1}Z})=0.69315$ and $H_{betP}(m_{{T1}Z})=H_{betP}(m_{{T2}Z})=0.69315$.\\
 
These first results indicate that  approaches based on absolute entropy-like measures appear to be useless for solving BAP since there is actually no reason which justifies that the correct assignment corresponds to the absolute minimum entropy-like measure just because $m_{Z}$ can stem from the least informational source. The association solution itself is actually independent of the informational content of each source.\\

An other attempt is to use rather the minimum of variation of entropy as decision criterion. Thus, the  following $\min\{\Delta_1(.),\Delta_2(.)\}$ criterions are examined; where variations $\Delta_i(.)$ for $i=1,2$ are defined as the
\begin{itemize}
\item variation of extended entropy:
$$\Delta_i(H_{ext})\triangleq H_{ext}(m_{{T_i}Z}) - H_{ext}(m_{{T_i}})$$
\item variation of generalized entropy:
$$\Delta_i(H_{gen}) \triangleq H_{gen}(m_{{T_i}Z}) - H_{gen}(m_{{T_i}})$$
\item variation of pignistic entropy:
$$\Delta_i(H_{betP}) \triangleq H_{betP}(m_{{T_i}Z}) -H_{betP}(m_{{T_i}})$$
\end{itemize}

Only the 2nd criterion, i.e. $\min(\Delta_i(H_{gen}))$ provides actually the correct solution 
for problem 1 and none of these criterions gives correct solution for problem 2.\\

The last idea is then to use the minimum of relative variations of pignistic probabilities of $\theta_1$ and $\theta_2$ given by the minimum on $i$ of 
$$\Delta_i(P) \triangleq \sum_{j=1}^2 \frac{|P_{T_iZ}(\theta_j)-P_{T_i}(\theta_j)|}{P_{T_i}(\theta_j)}$$
\noindent where $P_{T_iZ}(.)$ and $P_{T_i}(.)$ are respectively  the pignistic transformations of  $m_{T_iZ}(.)$ and $m_{T_i}(.)$.
Unfortunately, this criterion is unable to provide the solution for problems 1 and 2 because one has here in both problems $\Delta_1(P)=\Delta_2(P)=0$.

\subsection{The Schubert's approach}

We examine now  the possibility of using a Dempster-Shafer clustering method 
based on metaconflict function (MC-DSC) proposed in Johan Schubert's research works \cite{Schubert_1993,Schubert_2002} for solving the 
associations problems 1 and 2. A DSC method is a method of clustering uncertain data using the conflict in Dempster's rule as a distance measure  \cite{Schubert_2002a}. The basic idea is to separate/partition evidences by their conflict rather than by their proposition's event parts. Due to space limitation, we will just summarize here the principle of the {\it{classical}} MC- DSC method. 
Assume a given set of evidences (bba) $E(k)\triangleq\{m_{T_i}(.), i=1,\ldots,n\}$ is available at a given index (space or time or whatever) $k$ and suppose that a given set $E(k+1)\triangleq\{m_{z_j}(.), j=1,\ldots,m\}$ of new bba is then available for index $k+1$. The complete set of evidences representing all available information at index $k+1$ is $\chi=E(k) \cup E(k+1)\triangleq\{e_1,\ldots,e_q\}\equiv\{m_{T_i}(.), i=1,\ldots,n,m_{z_j}(.), j=1,\ldots,m\}$ with $q=n+m$. The problem we are faced now is to find the optimal partition/assignment of $\chi$ in disjoint subsets $\chi_p$ in order to combine informations within each $\chi_p$ in a coherent and efficient way.  The idea is to combine, in a first step, the set of bba belonging to the same subsets $\chi_p$ into a new bba $m_p(.)$ having a corresponding conflict factor $k_p$. The conflict factors $k_p$ are then used, in a second step, at a metalevel of evidence associated with the new frame of discernment $\Theta=\{AdP,\neg Adp\}$ where $AdP$ is short for {\it{adequate partition}}. From each subset $\chi_p$, $p=1,\ldots P$ of the partition under investigation, a new bba is defined as:
$$m_{\chi_p}(\neg AdP) \triangleq k_p \quad \text{and} \quad m_{\chi_p}(\Theta) \triangleq 1 - k_p$$
The combination of all these metalevel bba $m_{\chi_p}(.)$ by Dempster's rule yields a global bba 
$$m(.)=m_{\chi_1}(.)\oplus\ldots\oplus m_{\chi_P}(.)$$
\noindent with a corresponding {\sl{metaconflict factor}} denoted $Mcf(\chi_1,\ldots,\chi_P)\triangleq k_{1,\ldots,P}$. It can be shown \cite{Schubert_1993} that the metaconflict factor can be easily calculated directly from conflict factors $k_p$ by the following metaconflict function (MCF)
\begin{equation}
Mcf(\chi_1,\ldots,\chi_P)=1 - \prod_{p=1}^P (1-k_p)
\label{eqMCF}
\end{equation}

By minimizing the metaconflict function (i.e. by browsing all potential assignments), we intuitively expect to find the optimal/correct partition which will hopefully solve our association problem. Let's go back now to our very simple association problems 1 and 2 and examine the results obtained from the MC-DSC method.\\

The information available in association problems is denoted $\chi=\{m_{T_1}(.),m_{T_2}(.),m_Z(.)\}$.
We now examine all possible partitions of $\chi$ and the corresponding metaconflict factors and decision (based on minimum metaconflict function criterion) as follows:

\begin{itemize}
\item Analysis for  problem 1: 
\begin{itemize}
\item the (correct) partition $\chi_1=\{m_{T_1}(.),m_Z(.)\}$ and $\chi_2=\{m_{T_2}(.)\}$ yields
through Dempter's rule the conflict factors $k_1\triangleq k_{{T_1}Z}=0.5$ for subset $\chi_1$ and $k_2=0$ for subset $\chi_2$ since there is no combination at all (and therefore no conflict) in $\chi_2$. According to \eqref{eqMCF}, the value of the metaconflict is equal to
$$\text{Mcf}_1=1-(1-k_1)(1-k_2)=0.5 \equiv k_1$$
\item the (wrong) partition $\chi_1=\{m_{T_1}(.)\}$ and $\chi_2=\{m_{T_2}(.),m_Z(.)\}$ yields
the conflict factors $k_1=0$ for subset $\chi_1$ and $k_2=0.1$ for subset $\chi_2$. The value of the metaconflict is now equal to
$$\text{Mcf}_2=1-(1-k_1)(1-k_2)=0.1 \equiv k_2$$
\item since $\text{Mcf}_1>\text{Mcf}_2$, the minimum of the metaconflict function provides the wrong assignment and the MC-DSC approach fails to generate the solution for the problem 1.
\end{itemize}
\item Analysis for  problem 2:
\begin{itemize}
\item the (wrong) partition $\chi_1=\{m_{T_1}(.),m_Z(.)\}$ and $\chi_2=\{m_{T_2}(.)\}$ yields
through Dempter's rule the conflict factors $k_1\triangleq k_{{T_1}Z}=0.1$ for subset $\chi_1$ and $k_2=0$ for subset $\chi_2$ since there is no combination at all (and therefore no conflict) in $\chi_2$. According to \eqref{eqMCF}, the value of the metaconflict is equal to
$$\text{Mcf}_1=1-(1-k_1)(1-k_2)=0.1 \equiv k_1$$
\item the (correct) partition $\chi_1=\{m_{T_1}(.)\}$ and $\chi_2=\{m_{T_2}(.),m_Z(.)\}$ yields
the conflict factors $k_1=0$ for subset $\chi_1$ and $k_2=0.02$ for subset $\chi_2$. The value of the metaconflict is now equal to
$$\text{Mcf}_2=1-(1-k_1)(1-k_2)=0.02 \equiv k_2$$
\item since $\text{Mcf}_2<\text{Mcf}_1$, the minimum of the metaconflict function provides in this case the correct solution for the problem 2.
\end{itemize}
\end{itemize}

From these very simple examples, it is interesting to note that the Schubert's approach is actually exactly equivalent (in these cases) to the min-conflict approach detailed in section 3.2 and thus will not provide better results. It is also possible to show that the Schubert's approach also fails if one considers  
jointly the two observed bba $m_{Z_1}(.)$ and $m_{Z_2}(.)$ corresponding to problems 1 and 2 with $m_{T_1}(.)$ and $m_{T_2}(.)$. If one applies the principle of minimum metaconflict function, one will take the wrong decision since the wrong partition $\{(Z_1,T_2),(Z_2,T_1)\}$ will be declared. This result is in contradiction with our intuitive expectation for the true opposite partition $\{(Z_1,T_1),(Z_2,T_2)\}$ taking into account the coincidence of the respective belief functions. 

\section{A short DSmT presentation}
It has been reported in \cite{Dezert_2003,Lowrance_Garvey_1983,Schubert_1993} (and references therein) 
that the use of the DST must usually be done with extreme caution if one has to take a final and important 
decision from the result of the Dempter's rule of combination. In most 
of practical fusion applications based on the DST, some ad-hoc 
or heuristic techniques must always be added to the fusion process to 
manage or reduce the possibility of high degree of conflict between 
sources. Otherwise, the fusion results lead to a very dangerous 
conclusions (or cannot provide a reliable results at all). 
The practical limitations of the DST come essentially from its inherent
 following constraints which are closely related with the acceptance of the third exclude principle 
\begin{enumerate}
\item[(C1)] - the DST considers a discrete and finite frame of discernment 
$\Theta$ based on a set of exhaustive and exclusive elementary elements $\theta_i$.
\item[(C2)] - the bodies of evidence are assumed independent and provide their own belief function on the powerset $2^\Theta$ but with {\it{same interpretation}} for $\Theta$. 
\end{enumerate}

These two constraints therefore do not allow us to deal with the more general 
and practical problems involving uncertain reasoning and the fusion of 
uncertain, imprecise and paradoxical sources of information. To overcome these major limitations and drawbacks relative to the Dempster's rule of combination, a recent theory of plausible and paradoxical reasoning, called DSmT,  has been developed by Dezert and Smarandache in \cite{Dezert_2002b,Smarandache_2000,Smarandache_2002} and recently improved in \cite{Dezert_2003}.\\

The foundations of the DSmT is to refute the principle of the third exclude and to allow the possibility for paradoxes (partial vague overlapping) between elements of the frame of discernment.  The relaxation of the constraint C1 can be justified since the elements of $\Theta$ correspond generally only to imprecise/vague notions and concepts so that no refinement of $\Theta$ satisfying the first constraint is actually possible (specially if natural language is used to describe elements of $\Theta$). \\

The DSmT refutes also the excessive requirement imposed by C2 since it seems clear to us that, 
the {\sl{same}} frame $\Theta$ is usually interpreted differently by the distinct sources of evidence (experts). Some subjectivity on the information 
provided by a source of information is almost unavoidable, otherwise this would assume, as within the DST, that all bodies of 
evidence have an objective/universal (possibly uncertain) interpretation or measure of the phenomena under consideration which unfortunately rarely (never) occurs in reality, but when bba are based on some {\it{objective probabilities}} transformations (in such cases however probability theory tools become optimal tools to process all the available information; and the DST - as well as the DSmT - becomes useless). If we now get out of the probabilistic background argumentation, we claim that in most of cases, the sources of evidence provide their 
beliefs about some hypotheses only with respect to their own worlds of knowledge and experience without reference to the (inaccessible) absolute 
truth of the space of possibilities.\\

The DSmT  includes the possibility to deal with 
evidences arising from different sources of information which don't have 
access to absolute interpretation of the elements $\Theta$ under consideration and can be interpreted as a general 
and direct extension of probability theory and the Dempster-Shafer theory in the 
following sense. Let $\Theta=\{\theta_{1},\theta_{2}\}$ be the simplest frame of 
discernment involving only two elementary hypotheses (with no more assumptions on 
$\theta_{1}$ and $\theta_{2}$), then 
\begin{itemize}
\item Probability theory deals with probability assignments 
$m(.)\in [0,1]$ such that $m(\theta_{1})+m(\theta_{2})=1$
\item DST deals with bba $m(.)\in [0,1]$ such that 
$m(\theta_{1})+m(\theta_{2})+m(\theta_{1}\cup\theta_{2})=1$
\item DSmT theory deals with bba $m(.)\in [0,1]$ such that $m(\theta_{1})+m(\theta_{2})+m(\theta_{1}\cup\theta_{2})+m(\theta_{1}\cap\theta_{2})=1$
\end{itemize}

\subsection{Hyper-powerset and DSm rule}

Let $\Theta=\{\theta_{1},\ldots,\theta_{n}\}$ be a set of $n$ 
elements which cannot be precisely defined and separated so that no 
refinement of $\Theta$ in a new larger set $\Theta_{ref}$ of disjoint elementary 
hypotheses is possible. The {\sl{hyper-power}} set $D^\Theta$ is defined as the set of all 
composite possibilities build from $\Theta$ with $\cup$ and $\cap$ 
operators such that $\forall A\in D^\Theta, B\in D^\Theta, (A\cup 
B)\in D^\Theta$ and $(A\cap B)\in D^\Theta$. The cardinality of $D^\Theta$ is majored by 
$2^{2^n}$ when $\text{Card}(\Theta)=\mid\Theta\mid =n$. The generation 
of hyper-power set $D^\Theta$ is closely related with the famous Dedekind's problem on enumerating the 
set of monotone Boolean functions. An algorithm for generating $D^\Theta$ can 
be found in \cite{Dezert_Smarandache_2003} for convenience. From a general frame of discernement $\Theta$, we define a map $m(.): 
D^\Theta \rightarrow [0,1]$ associated to a given source of evidence $\mathcal{B}$ 
which can support paradoxical information, as follows
\begin{equation*}
m(\emptyset)=0 \qquad \text{and}\qquad \sum_{A\in D^\Theta} m(A) = 1 
\end{equation*}
The quantity $m(A)$ is called $A$'s {\sl{general basic belief 
number}} (gbba) or the general basic belief mass for $A$.
The belief and plausibility functions are defined in almost the same manner as within the DST, i.e.
\begin{equation*}
\text{Bel}(A) = \sum_{B\in D^\Theta, B\subseteq A} m(B)
\end{equation*}
\begin{equation*}
\text{Pl}(A) = \sum_{B\in D^\Theta, B\cap A\neq\emptyset} m(B)
\end{equation*}
Note that the {\it{classical}} complementary $A^c$ of any given proposition $A$ is not involved within DSmT just because of the refutation of the third exclude principle.
These definitions are compatible with the DST definitions when the sources of 
information become uncertain but rational (they do not support paradoxical 
information). We still have $\forall A\in D^\Theta, \text{Bel}(A)\leq \text{Pl}(A)$.\\

The DSm rule of combination $m(.)\triangleq [m_{1}\oplus m_{2}](.)$ of two distinct (but potentially paradoxical) sources of evidences $\mathcal{B}_{1}$ and  $\mathcal{B}_{2}$ over the same 
general frame of discernment $\Theta$ with belief functions $\text{Bel}_{1}(.)$ and 
 $\text{Bel}_{2}(.)$ associated with general information granules $m_{1}(.)$ and $m_{2}(.)$ is  then given by $\forall C\in D^\Theta$,
 \begin{equation*}
m(C) = 
 \sum_{A,B\in D^\Theta, A\cap B=C}m_{1}(A)m_{2}(B)
 \label{JDZT}
 \end{equation*}
Since $D^\Theta$ is closed under $\cup$ and $\cap$ operators, this new rule 
of combination guarantees that $m(.): D^\Theta \rightarrow [0,1]$ is a proper general information granule. 
This rule of combination is commutative and associative 
and can always be used for the fusion of paradoxical or rational sources of 
information (bodies of evidence). It is important to note that any fusion of sources of information 
generates either uncertainties, paradoxes or {\sl{more generally both}}. This is intrinsic to the 
general fusion process itself. The theoretical justification of the DSm rule can be found in \cite{Dezert_2003}. As within the DST framework, it is possible to build a subjective probability 
 measure $P^\star\{.\}$ from the bba $m(.)$ with  the generalized pignistic transformation (GPT) \cite{Dezert_2003,Dezert_Smarandache_Daniel_2003} defined $\forall A \in D^\Theta$ by,
\begin{equation*}
P^\star\{A\}=\sum_{C \in D^\Theta | A\cap C\neq \emptyset}  \frac{\mathcal{C}_{\mathcal{M}^f}(C\cap A)}{\mathcal{C}_{\mathcal{M}^f}(C)}m(C)
\end{equation*}
where $\mathcal{C}_{\mathcal{M}^f}(X)$ denotes the DSm cardinal of proposition $X$ for the free-DSm model $\mathcal{M}^f$ of the problem under consideration here \cite{Dezert_2003f}.
From any generalized bba $m(.)$ and its corresponding pignistic transformation $P^\star(.)$, one can also define the following new entropy measures
\begin{itemize}
\item New extended entropy-like measure:
$$H^\star_{ext}(m)\triangleq - \sum_{A\in D^\Theta} m(A)\log(m(A))$$
\item New generalized pignistic entropy :
\begin{equation*}
H_{betP}^\star(P^\star)\triangleq  -\sum_{A\in \mathcal{V}} P^\star\{A\}\ln(P^\star\{A\})
\end{equation*}
\end{itemize}
\noindent
where $\mathcal{V}$ denotes the parts of the Venn diagram of the model $\mathcal{M}^f$.

\section{DSmT approaches for BAP}

As within DST, several approaches can be attempted to try to solve the Blackman's Association problems (BAP). The first attempts are based on the minimum on $i$ of new extended entropy-like measures $H^\star_{ext}(m_{{T_i}Z})$ or on the minimum $H^\star_{betP}(P^\star)$. Both approaches actually fail for the same reason as for the DST-based minimum entropy criterions.\\

The second attempt is based on the minimum of variation of the new entropy-like measures as criterion for the choice of the decision with  the new extended entropy-like measure:
$$\Delta_i(H^\star_{ext})\triangleq H^\star_{ext}(m_{{T_i}Z}) - H^\star_{ext}(m_{{T_i}})$$
or the new generalized pignistic entropy:
$$\Delta_i(H^\star_{betP}) \triangleq H^\star_{betP}(P^\star\{.\vert m_{{T_i}Z}\}) - H^\star_{betP}(P^\star\{.\vert m_{T_i} \})$$

The min. of $\Delta_i(H^\star_{ext})$ gives us the wrong solution for problem 1 since $\Delta_1(H^\star_{ext})=0.34657$ and
$\Delta_2(H^\star_{ext})=0.30988$ while min. of $\Delta_i(H^\star_{betP})$ give us the correct solution since $\Delta_1(H^\star_{betP})=-0.3040$ and $\Delta_2(H^\star_{betP})=-0.0960$. Unfortunately, both the $\Delta_i(H^\star_{ext})$ and $\Delta_i(H^\star_{betP})$ criterions fail to provide the correct solution for problem 2 since one gets $\Delta_1(H^\star_{ext})=0.25577 < \Delta_2(H^\star_{ext})=0.3273$ and $\Delta_1(H^\star_{betP})=-0.0396 < \Delta_2(H^\star_{betP})=-0.00823$.\\

The third proposed approach is to use the criterion of the minimum of relative variations of pignistic probabilities of $\theta_1$ and $\theta_2$ given by the minimum on $i$ of 
$$\Delta_i(P^\star) \triangleq \sum_{j=1}^2 \frac{|P^\star_{T_iZ}(\theta_j)-P^\star_{T_i}(\theta_j)|}{P^\star_{T_i}(\theta_j)}$$

This third approach fails to find the correct solution for problem 1 
(since $\Delta_1(P^\star)=0.333 > \Delta_2(P^\star)=0.268$) but succeeds to get the correct solution for problem 2
(since $\Delta_2(P^\star)=0.053 < \Delta_1(P^\star)=0.066$).\\

The last proposed approach is based on relative variations of  pignistic probabilities conditioned
by the correct assignment. The criteria is defined as the minimum of 
$$\delta_i(P^\star) \triangleq \frac{|\Delta_i(P^\star | Z) - \Delta_i(P^\star | \hat{Z}=T_i)|}{\Delta_i(P^\star | \hat{Z}=T_i)}$$
\noindent
where  $\Delta_i(P^\star | \hat{Z}=T_i)$ is obtained as for $\Delta_i(P^\star)$ but by forcing $Z=T_i$ or equivalently $m_Z(.)=m_{T_i}(.)$ for the derivation of pignistic probabilities $P^\star_{T_iZ}(\theta_j)$. This last criterion yields the correct solution for problem 1 (since $\delta_1(P^\star)=|0.333-0.333|/0.333=0 < \delta_2(P^\star)=|0.268-0.053|/0.053\approx 4$) and simultaneously for problem 2 (since $ \delta_2(P^\star)=|0.053-0.053|/0.053=0 < \delta_1(P^\star)=|0.066-0.333|/0.333\approx 0.8$).
 
\section{Monte-Carlo simulations}

As shown on the two previous BAP, it is difficult to find a general method for solving both these particular (noise-free $m_Z$) BAP and all general problems involving noisy attribute bba $m_Z(.)$. The proposed methods have been examined only for the original BAP and no {\sl{general}} conclusion can be drawn from our previous analysis about the most efficient approach.
The evaluation of the global performances/efficiency of previous approaches can however be estimated quite easily through 
Monte-Carlo simulations. Our Monte-carlo simulations are based on 50.000 independent runs and have been done both for the noise-free case (where $m_Z(.)$ matches perfectly with either $m_{T_1}(.)$ or $m_{T_2}(.)$) and for two noisy cases (where $m_Z(.)$ doesn't match perfectly one of the predicted bba). Two noise levels (low and medium) have been tested for the noisy cases. A basic run consists in generating randomly the two predicted bba $m_{T_1}(.)$ and $m_{T_2}(.)$ and an observed bba $m_Z(.)$ according  to a random assignment $m_Z(.) \leftrightarrow m_{T_1}(.)$ or $m_Z(.) \leftrightarrow m_{T_2}(.)$. Then we evaluate the percentage of right assignments for all chosen association criterions described in this paper. The introduction of  noise on perfect (noise-free) observation $m_Z(.)$ has been obtained by the following procedure (with notation $A_1\triangleq \theta_1$, $A_2\triangleq \theta_2$ and
$A_2\triangleq \theta_1\cup\theta_2$): $m_Z^{\text{\tiny noisy}}(A_i)=\alpha_i m_Z(A_i)/K$ where $K$ is a normalization constant such as $\sum_{i=1}^3 m_Z^{\text{\tiny noisy}}(A_i)=1$ and weighting coefficients $\alpha_i \in [0;1]$ are given by $\alpha_i=1/3 \pm \epsilon_i$ such that  $\sum_{i=1}^3 \alpha_i=1$. \\

The table 1  shows the Monte-Carlo results obtained with all investigated criterions for the following 3 cases: noise-free (NF), low noise (LN) and medium noise (MN) related to the observed bba $m_Z(.)$. The two first rows of the table correspond to simplest approach. The next twelve rows correspond to DST-based approaches.
\begin{center}
\begin{tabular}{|l|r|r|r|}
\hline
Assoc. Criterion &  NF & LN & MN \\
\hline
Min $d_{L^1}(T_i,Z)$ & 100 & 97.98 & 92.14\\
Min $d_{L^2}(T_i,Z)$ & 100 & 97.90 & 92.03\\
\hline
Min $k_{{T_i}Z}$                    & 70.01 & 69.43 & 68.77\\
Min $L(Z|T_i)$                        & 70.09 &  69.87 & 67.86\\
Min $d_{L^1}(T_i,T_iZ)$      & 57.10 & 57.41 & 56.30\\
Min $d_{L^2}(T_i,T_iZ)$      & 56.40 & 56.80 & 55.75\\
Min $H_{ext}(m_{{T_i}Z})$   & 61.39 & 61.68 & 60.85 \\
Min $H_{gen}(m_{{T_i}Z})$ & 58.37 & 58.79 & 57.95\\
Min $H_{betP}(m_{{T_i}Z})$ & 61.35 & 61.32 & 60.34\\
Min $\Delta_i(H_{ext})$ & 57.66 & 56.97 & 55.90\\
Min $\Delta_i(H_{gen})$ & 57.40 & 56.80 & 55.72\\
Min $\Delta_i(H_{betP})$ & 71.04 & 69.15 & 66.48\\
Min $\Delta_i(P)$ & 69.25 & 68.99 & 67.35\\
Min $Mcf_i$          & 70.1 & 69.43& 68.77\\
\hline
\end{tabular}
\end{center}
\begin{center}
Table 1 : $\%$ of success of association methods
\end{center}

The table 2  shows the Monte-Carlo results obtained for the 3 cases: noise-free (NF), low noise (LN) and medium noise (MN) related to the observed bba $m_Z(.)$ with the DSmT-based approaches.

\begin{center}
\begin{tabular}{|l|r|r|r|}
\hline
Assoc. Criterion &  NF & LN & MN \\
\hline
\hline
Min $H^\star_{ext}(m_{{T_i}Z})$ & 61.91 & 61.92 & 60.79\\
Min $H^\star_{betP}(P^\star)$ & 42.31 & 42.37 & 42.96\\
Min $\Delta_i(H^\star_{ext})$ & 67.99 & 67.09 & 65.72\\
Min $\Delta_i(H^\star_{betP})$ & 42.08 & 42.11 &42.21\\
Min $\Delta_i(P^\star)$ & 76.13 & 75.3 & 72.80\\
Min $\delta_i(P^\star)$ & 100 & 90.02& 81.31\\
\hline
\end{tabular}
\end{center}
\begin{center}
Table 2 : $\%$ of success of DSmT-based methods
\end{center}

\section{Conclusion}

A deep examination of the Blackman's association problem has been presented. Several methods have been proposed and compared through Monte Carlo simulations. Our results indicate that the commonly used min-conflict method doesn't provide the best performance in general (specially w.r.t. the simplest distance approach). Thus the metaconflict approach, equivalent here to min-conflict, does not allow to get the optimal efficiency. The Blackman's approach and min-conflict give same performances. All entropy-based methods are less efficient than the min-conflict approach. More interesting, from  the results based on the generalized pignistic entropy approach, the entropy-based methods seem actually not appropriate for solving BAP since there is no fundamental reason to justify them. The min-distance approach of Tchamova is the least efficient method among all methods when abandoning entropy-based methods. Monte carlo simulations have shown that only methods based on the relative variations of generalized pignistic probabilities build from the DSmT outperform all methods examined in this work but the simplest one.

\end{document}